\documentclass[12pt, 14paper,reqno]{amsart}
\usepackage{amsmath,amsfonts,amssymb}
\usepackage[breaklinks]{hyperref}
\usepackage{graphicx}
\usepackage{longtable}
\makeatletter
\@namedef{subjclassname@2020}{%
	\textup{2020} Mathematics Subject Classification}
\makeatother
\theoremstyle{plain}
\newtheorem{thm}{Theorem}[section]
\newtheorem{lem}{Lemma}[section]%[section]https://www.overleaf.com/project/656efb8c0bb2b89cedb1fe05

\newtheorem{q}{Question}[section]
\newtheorem{conj}{Conjecture}[section]

\newtheorem{thma}{Theorem}

\theoremstyle{proof}

\numberwithin{equation}{section}

\begin{document} 
	\title[$D(n)$-quadruples in $\mathbb{Z}(\sqrt{4k+2})$]{Diophantine $D(n)$-quadruples in $\mathbb{Z}[\sqrt{4k + 2}]$}
	\author{Kalyan Chakraborty, Shubham Gupta and Azizul Hoque}
	\address{KC @Department of Mathematics, SRM University AP, Neerukonda, Mangalagiri, Guntur-522240, Andhra Pradesh, India.}
	\email{kalychak@gmail.com}
	\address{SG @Harish-Chandra Research Institute,  A CI of Homi Bhabha National
		Institute, Chhatnag Road, Jhunsi, Prayagraj - 211019, India.}
	\email{shubhamgupta2587@gmail.com}
	\address{AH @Department of Mathematics, Faculty of Science, Rangapara College, Rangapara, Sonitpur-784505, Assam, India.}
	\email{ahoque.ms@gmail.com}
	\keywords{Diophantine quadruples; Pellian equations; Quadratic fields}
	\subjclass[2020] {11D09; 11R11}
	\date{\today}
	\maketitle
	
	\begin{abstract}
		
		Let $d$ be a square-free integer and $\mathbb{Z}[\sqrt{d}]$  a quadratic ring of integers. For a given $n\in\mathbb{Z}[\sqrt{d}]$,  a set of $m$ non-zero distinct elements in $\mathbb{Z}[\sqrt{d}]$ is called a Diophantine $D(n)$-$m$-tuple (or simply $D(n)$-$m$-tuple) in $\mathbb{Z}[\sqrt{d}]$ if product of any  two of them plus $n$ is a square in $\mathbb{Z}[\sqrt{d}]$. Assume that $d \equiv 2 \pmod 4$ is a positive integer such that $x^2 - dy^2 = -1$ and $x^2 - dy^2 = 6$ are solvable in integers. In this paper, we prove the existence of  infinitely many $D(n)$-quadruples in $\mathbb{Z}[\sqrt{d}]$ for $n = 4m + 4k\sqrt{d}$ with $m, k \in  \mathbb{Z}$ satisfying $m \not\equiv 5 \pmod{6}$ and $k \not\equiv 3 \pmod{6}$. Moreover, we prove the same for $n = (4m + 2) + 4k\sqrt{d}$ when either $m \not\equiv 9 \pmod{12}$ and $k \not\equiv 3 \pmod{6}$, or $m \not\equiv 0 \pmod{12}$ and $k \not\equiv 0 \pmod{6}$. At the end, some examples supporting the existence of quadruples in $\mathbb{Z}[\sqrt{d}]$ with the property $D(n)$ for  the above exceptional $n$'s are provided for $d = 10$.
		
	\end{abstract}
	
	\section{Introduction}
	A set $\{a_1, a_2, \ldots, a_m\}$ of $m$ distinct positive integers is called a Diophantine $m$-tuple with the property $D(n)$ (or simply $D(n)$-$m$-tuple) for a given non-zero integer $n$, if  $a_ia_j+n$ is a perfect square for all $1\leq i<j\leq m$. For $n = 1$, such an $m$-tuple is called Diophantine $m$-tuple instead of Diophantine $m$-tuple with the property $D(1)$. 
	The question of constructing such tuples was first studied by Diophantus of Alexandria, who found a Diophantine  quadruple of rationals  $\{1/16, 33/16, 17/4, 105/16\}$ with the property $D(1)$.
 However, it was Fermat who first found a Diophantine quadruple $\{1, 3, 8, 120\}$ in integers. Later, Baker and Davenport \cite{BD1969} proved that Fermat's quadruple can not be extended to Diophantine quintuple. Dujella \cite{DU2004} proved the non-existence of Diophantine sextuple and that there are at most finitely many integer Diophantine quintuples. Recently, He, Togb\'e and Ziegler \cite{HTZ2019} proved the non-existence of integer Diophantine quintuples, and in this way, they  solved a long-standing open problem. 
	On the other hand, Bonciocat, Cipu and Mignotte \cite{BCM2020} proved a conjecture of Dujella \cite{DU1993}, which states that there are no $D(-1)$-quadruples. It is also known due to  Trebje\v{s}anin and Filipin \cite{BTF2019}  that there do not exist  $D(4)$-quintuples. A brief survey on this topic can be found in \cite{DU23}. We also refer \cite{BR1985, CGH22, DU21, DU2024, EFF2014} to the reader for more information about $D(n)$-$m$-tuples.
	 
 Let $\mathcal{R}$ be a commutative ring with unity. For a given $n\in\mathcal{R}$,  a set $\{a_1, a_2, \ldots, a_m\} \subset \mathcal{R} \setminus \{0\}$ is called a Diophantine $m$-tuple  with the property $D(n)$ in $\mathcal{R}$ (or simply $D(n)$-$m$-tuple in $\mathcal{R}$), if  $a_i a_j+n$ is a perfect square in $\mathcal{R}$ for  all $1 \leq i < j \leq m$. Let $K$ be an imaginary quadratic number field and $\mathcal{O}_K$ be its ring of integers. In 2019,  Ad\v{z}aga \cite{A2019} proved that there are no $D(1)$-$m$-tuples  in $\mathcal{O}_K$ when $m\geq 42$. Recently, Gupta \cite{G2021} proved that there do not exist $D(-1)$-$m$-tuple for $m\geq 37$.
 %\pagebreak\\
 It is interesting to note that $D(n)$-quadruples are related to the representations of $n$ by the binary quadratic form $x^2-y^2$.
	In particular, Dujella \cite{DU1993} proved  that a  $D(n)$-quadruple in integers exists if and only if $n$ can be written as a difference of two squares, up to finitely many exceptions. Later, Dujella \cite{DU1997} proved the above fact in Gaussian integers. Further, the above fact also holds for the ring of integers of  $\mathbb{Q}(\sqrt{d})$ for certain $d \in \mathbb{Z}$ (see, \cite{FR2004, FR2008, FR2009, FS2014, MR2004,  SO2013}). These results motivated Franu\v{s}i\'c and Jadrijevi\'c to post the following conjecture:
\begin{conj}[{\cite[Conjecture 1]{FJ2019}}]\label{Con1.1}
		Let $\mathcal{R}$ be a commutative ring with unity $1$ and $n \in \mathcal{R}\setminus \{0\}$. Then a $D(n)$-quadruple  exists if and only if  $n$ can be written as a difference of two squares in $\mathcal{R}$, up to finitely many exceptions of $n$.
\end{conj}
This conjecture was verified for rings of integers of certain number fields (cf. \cite{FR2004, FR2008, FR2009, FR2013, FJ2019, FS2014, MA2012, MR2004,  SO2013}).

 The following  notations will be followed throughout the paper. 
	
	\begin{itemize}
		\item $(a, b)= a + b\sqrt{d}$,
		\item  $k(a, b) = (ka, kb)$ for $k\in \mathbb{Z}$,
		\item 
		Let $\alpha=(a, b)$. The norm $\text{Nm}$ of $\alpha$ is given by 
		$$\text{Nm}(\alpha) := (a, b)(a, -b),$$
		\item $(x, y) \equiv (a, b) \text{~(mod~} (c, e))$ means that
		$
		x \equiv a \pmod{c} \text{~~~and~~~} y \equiv b \pmod{e}.
		$
		
	\end{itemize}
In the rest of paper, we fix $d \equiv 2 \pmod{4}$ to be a square-free positive integer. We set $\mathcal{S}$ and $\mathcal{T}$ in $\mathbb{Z}[\sqrt{d}]$ as follows:
	\begin{align*}
		\mathcal{S} :=& \{(4m, 4k + 1), (4m, 4k + 2), (4m, 4k + 3), (4m + 1, 4k + 1), (4m + 1, 4k + 3), (4m + 2, \\ & 4k + 1), (4m + 2, 4k + 3), (4m + 3, 4k + 1),(4m + 3, 4k + 3)\},\\ 
		\mathcal{T} :=& \{(4m, 4k), (4m + 1, 4k), (4m + 1, 4k + 2), (4m + 2, 4k), (4m + 2, 4k + 2), (4m + 3, 4k),\\ & (4m + 3, 4k + 2)\},
	\end{align*}
	where $m, k \in \mathbb{Z}$.
It is easy to check that if $n \in \mathbb{Z}[\sqrt{d}]$ then $n \in \mathcal{S} \cup \mathcal{T}$. In \cite{FR2004}, Franu\v si\'c  proved that there does not exist any $D(n)$-quadruple in $\mathbb{Z}[\sqrt{d}]$ for $n \in \mathcal{S}$.

Thus, it is natural to ask `{\it whether there exists any Diophantine quadruple in $\mathbb{Z}[\sqrt{d}]$ for $n\in \mathcal{T}$}'. 
Very recently, in \cite{CGH2022} the present authors answered this question for $n \in \mathcal{T}\setminus \{(4m, 4k), (4m + 2, 4k)\}$. More precisely, the authors proved the following result: 

\begin{thma}[{\cite[Theorem 1.1]{CGH2022}}]\label{thmA}
Assume that $d\equiv 2\pmod 4$ is a square-free positive integer and the equations \eqref{eqi1} and \eqref{eqi2} are solvable. Then there exist infinity many quadruples in $\mathbb{Z}[\sqrt{d}]$ with the property $D(n)$ when $n \in \{(4m + 1) + 4k\sqrt{d}, (4m + 1) + (4k + 2)\sqrt{d}, (4m + 3) + 4k\sqrt{d}, (4m + 3) + (4k + 2)\sqrt{d}, (4m + 2) + (4k + 2)\sqrt{d}\}$ with $m, k \in \mathbb{Z}$. 
\end{thma}
As a consequence of Theorem \ref{thmA}, we were able to construct some counter examples of Conjecture \ref{Con1.1}. Namely, if $d=10$ and $n = 26 + 6\sqrt{10}$  or $d=58$ and $n = 18 + 2\sqrt{58}$, one can easily see that $n$ can not be represented as a difference of two squares in  $\mathbb{Z}[\sqrt{d}]$, but there exists a $D(n)$-quadruple in $\mathbb{Z}[\sqrt{d}]$. 

 In this paper, we consider the above mentioned problem for the remaining values of $n$. 
Let $d \equiv 2 \pmod 4$ be a square-free positive integer such that
	\begin{equation}\label{eqi1}
		x^2 - dy^2 = -  1
	\end{equation}
	and 
	\begin{equation}\label{eqi2}
		x^2 - dy^2 =   6
	\end{equation}
are solvable in integers. We prove the following results:

	\begin{thm}\label{1.2}
Let $d\equiv 2\pmod 4$ be a square-free positive integer such that \eqref{eqi1} and \eqref{eqi2} are solvable in integers. Let $n = (4m, 4k)$ with $m, k \in \mathbb{Z}$ such that $(m, k) \not\equiv (5, 3) \pmod {(6, 6)}$.  Then there exist infinitely many $D(n)$-quadruples in $\mathbb{Z}[\sqrt{d}]$.
\end{thm}

\begin{thm}\label{1.3} Let $d$ be as in Theorem \ref{1.2}. Then for $n = (4m + 2, 4k)$ with $m, k\in \mathbb{Z}$,  
there exist infinitely many $D(n)$-quadruples  in $\mathbb{Z}[\sqrt{d}]$ such that $(m, k) \not\equiv (9, 3),  (0, 0) \pmod {(12, 6)}$.
\end{thm}

In 1996,  Dujella \cite{DUJE1996} obtained several two-parameter polynomial families for quadruples with the property $D(n)$. Our proofs use the technique presented in \cite{DUJE1996}.

\section{Preliminaries}

We begin this section with the following lemma that follows from the definition of $D(n)$-quadruples in $\mathbb{Z}[\sqrt{d}]$.
	\begin{lem}\label{lemp}
Let  $\{a_1, a_2, a_3, a_4\}$ be a $D(n)$-quadruple. Then for any non-zero $w\in \mathbb{Z}[\sqrt{d}]$, with a square-free integer $d$, the set  $\{wa_1, wa_2, wa_3, wa_4\}$ is a  $D(w^2n)$-quadruple in $\mathbb{Z}[\sqrt{d}]$.
\end{lem}
The next lemma  helps us to  find the conditions under which the set  $\{a, b, a + b + 2r, a + 4b + 4r\}$ forms a $D(n)$-quadruple  in $\mathbb{Z}[\sqrt{d}]$ for any $n \in \mathbb{Z}[\sqrt{d}]$. 
	
\begin{lem}[{\cite[Lemma 2.5]{CGH2022}}]\label{lem2.1}
 The set $\{a, b, a + b + 2r, a + 4b + 4r\}$ of non-zero and distinct elements is a   $D(n)$-quadruple in $\mathbb{Z}[\sqrt{d}]$ for any $n\in\mathbb{Z}[\sqrt{d}]$,  if
$ ab + n = r^2$ and $3n = \alpha_1\alpha_2$ with $\alpha_1=a + 2r + \alpha$ and $ \alpha_2= a + 2r - \alpha$,
for some $ a, b, r, \alpha\in \mathbb{Z}[\sqrt{d}]$.
\end{lem}
The next two lemmas help us to apply Lemma \ref{lem2.1} in the proofs of Theorems \ref{1.2} and \ref{1.3}. Lemma \ref{lem3.1} is useful for the factorization of $3n$ in $\mathbb{Z}[\sqrt{d}]$, while Lemma \ref{lemm1} is useful to verify that the elements thus found are distinct and non-zero.
	
	\begin{lem}[{\cite[Lemma 3.1]{CGH2022}}]\label{lem3.1} Let $d\equiv 2\pmod 4$ be a square-free integer such that \eqref{eqi1} and \eqref{eqi2} are solvable in integers.    
		Then in $\mathbb{Z}[\sqrt{d}]$, the following statements hold:
		\begin{itemize}
			\item[(i)] elements of norm  $1$ have the form $(6a_1 \pm 1, 6b_1)$ and there are infinitely many of them;
			\item[(ii)] elements of norm  $-1$ have the form $(6a_1 \pm 3, 6b_1 \pm 1)$ and there are infinitely many such elements;
			\item[(iii)] $d \equiv 10 \pmod{48}$;
			\item[(iv)] elements of norm  $6$ have the form $(12M \pm 4, 6N \pm 1)$ and there are infinitely many such elements;
			\item[(v)] elements of norm  $-6$ have the form $(12M \pm 2, 6N \pm 1)$ and there are infinitely many such elements;
		\end{itemize}
		where $a_1, b_1, M ~\mbox{and} ~N \in \mathbb{Z}$.
	\end{lem}

	\begin{lem}[{\cite[Lemma 2.4]{CGH2022}}]\label{lemm1} Assume that $a_1, a_2, b_1, b_2, c_1, c_2, d_1, d_2, e_1 \in \mathbb{Z}$ with $a_1, a_2, b_1 \neq 0$. Then
		the following system of simultaneous equations
		\begin{align}\label{eQ7}
			\begin{cases}
				a_1x^2 + b_1y^2 + c_1x + d_1y + e_1 = 0,\\
				a_2xy + b_2x + c_2y + d_2 = 0
			\end{cases}
		\end{align}
		has only finitely many solutions in integers.
	\end{lem} 
	
	\section{Proof of Theorem \ref{1.2}}\label{s3}
	We first factorize $3n$ by using Lemmas \ref{lem2.1} and \ref{lem3.1}. We then use this factorization together with Lemma \ref{lem2.1} to construct Diophantine quadruples of certain forms with the property $D(n)$ under the condition of non-zero and distinctness. Finally these conditions are verified by using Lemma \ref{lemm1}.

	Here, $n = (4m, 4k)$ with $m, k\in \mathbb{Z}$. Thus
	$3n = 3(4m, 4k) = 6(2m, 2k)$
	and we choose $\alpha_1 = 6$ and $\alpha_2 = (2m, 2k)$ ($\alpha_1$ and $\alpha_2$ as in Lemma \ref{lem2.1}). Now Lemma \ref{lem2.1} entails, 
	\begin{equation}\label{eq5.1}
		a + 2r = (m + 3, k).
	\end{equation}
	We divide the proof into four cases based on the parity of $m$ and $k$.
	\subsection*{Case I:  Both $m$ and $k$ are even}
	Let $a = (6a_1 + 1, 6b_1)$ with $a_1, b_1\in \mathbb{Z}$ such that $\text{Nm}(a)=1$. Then by (i) of Lemma \ref{lem3.1}, there exist infinitely many  such $a$'s, and  \eqref{eq5.1} can be written as
	$$
	r = (m/2 + 1 - 3a_1, k/2 - 3b_1).
	$$
	As both  $m$ and $k$ are even, so $r \in \mathbb{Z}[\sqrt{d}]$. We employ these $a$ and $r$ in the equation $ab + n = r^2$ (as in Lemma \ref{lem2.1}) to get:
	$$
	b = ((m/2 + 1 - 3a_1)^2 + d(k/2 - 3b_1)^2 - 4m, 2(m/2 + 1 - 3a_1)(k/2 - 3b_1) - 4k)(6a_1 + 1, -6b_1).
	$$ 
	These choices of $a, b$ and $r$ give us infinitely many $D(n)$-quadruples $\{a,b,a+b+2r, a+4b+4r\}$ in $\mathbb{Z}[\sqrt{d}]$. Non-zero  and distinctness of these elements can easily be verified by Lemma \ref{lemm1}. 
	
	\subsection*{Case II: $m$ is odd and $k$ is even} As in Case I, we choose $a = 2(6a_1 + 1, 6b_1)$ with $a_1, b_1\in \mathbb{Z}$ and $\text{Nm}(a) = 4$. Then \eqref{eq5.1} gives,
	$$
	2r = (m + 1 - 12a_1, k - 12b_1).
	$$
	We write $m = 2m_1 + 1$ and $k = 2k_1$ for some $m_1, k_1 \in \mathbb{Z}$.  Then
	$$
	r = (m_1 + 1 - 6a_1, k_1 - 6b_1),
	$$
	which gives 
	$$
	b = \dfrac{1}{2}\left((m_1 + 1 - 6a_1)^2 + d(k_1 - 6b_1)^2 - 4m, 2(m_1 + 1 - 6a_1)(k_1 - 6b_1) - 4k\right)(6a_1 + 1, -6b_1).
	$$
We are looking for $b$ satisfying $b \in \mathbb{Z}[\sqrt{d}]$, so that $m_1$ should be odd and  $k_1$ should be even. These choices of $a, b$ and $r$ provide infinitely many $D(n)$-quadruples of the form $\{a,b,a+b+2r, a+4b+4r\}$ in $\mathbb{Z}[\sqrt{d}]$.
	
On the other hand for even $m_1$, we choose $a = 4(6a_1 + 1, 6b_1)$ with $a_1, b_1\in \mathbb{Z}$ and $\text{Nm}(a)=16$. Then as before we get 
	$$
	r = (m_1 - 12a_1, k_1 - 12b_1),
	$$
which provides 
	$$
	b = \dfrac{1}{4}\left((m_1 - 12a_1)^2 + d(k_1 - 12b_1)^2 - 4m, 2(m_1 - 12a_1)(k_1 - 12b_1) - 4k\right)(6a_1 + 1, -6b_1).
	$$
	Clearly $b \in \mathbb{Z}[\sqrt{d}]$ when $k_1$ is even. These give the required elements $a$, $b$ and $r$. Utilizing Lemma \ref{lem2.1}, this implies that the set $\mathcal{A} = \{a,b,a+b+2r, a+4b+4r\}$ forms a Diophantine quadruple in $\mathbb{Z}[\sqrt{d}]$ with the property $D(n)$, under the condition that all the elements of $\mathcal{A}$ must be non-zero and distinct from each other. These conditions can be verified by using Lemma \ref{lemm1}, except $a + 4b + 4r \ne 0$ and $a + 2r \ne 0$. We handle these exceptions separately since they do not fit into Lemma \ref{lemm1}.  We first consider $a + 2r = 0$. This gives $m_1 = -2$ and $k_1 = 0$. This gives $n  = -12$. Now if $a + 4b + 4r = 0$, then $(m_1, k_1) = (0, 0)$ or $(m_1, k_1) = (4, 0)$. This gives $n = 1, 36$, which are already known. 
	
The case $n = -12$ gives $
	3n = -18 \times 2
	$. We now choose  $\alpha_1 = -18$ and $\alpha_2 = 2$. 
As before, we  choose $a = 4(6a_1 + 1, 6b_1)$ with $a_1, b_1\in\mathbb{Z}$ and $\text{Nm}(a) = 16$, and thus
	$
	r = (-2 - 12a_1, -12b_1).
	$
	This gives $$
	b = ((1 + 6a_1, 6b_1)^2 + 3)(6a_1 + 1, -6b_1).
	$$
	Owing to the guaranteed existence  of infinitely many $a$'s,  there exist infinitely many $D(n)$-quadruples. 
	
	The  possibility of $m_1$ even and $k_1$ odd needs to be examined. In this case $n = (16m + 4, 16k + 8) = 2^2(4m + 1, 4k + 2)$, and thus the existence of infinitely many $D(n)$-quadruples in $\mathbb{Z}[\sqrt{d}]$ is guaranteed by \cite[Theorem 1.1]{CGH2022} and  Lemma \ref{lemp}. 
	
	\subsection*{Case III: $m$ is even and $k$ is odd}
	In this case, we consider
	$
	a = (6a_1 + 3, 6b_1 + 1)
	$ 
	with $a_1, b_1\in\mathbb{Z}$ and $\text{Nm}(a)=-1$. This provides us 
	$$
	b = ((m/2 - 3a_1)^2 + d((k - 1)/2 - 3b_1)^2 - 4m, 2(m/2 - 3a_1)((k - 1)/2 - 3b_1) - 4k)(-6a_1 - 3, 6b_1 + 1),
	$$
	(for the value of $r$ we use \eqref{eq5.1}). 
	As dealt with in the previous cases, these values of $a, b, r$ will guarantee infinitely many $D(n)$- quadruples in $\mathbb{Z}[\sqrt{d}]$.
	
	\subsection*{Case IV: Both $m$ and $k$ are odd}
	This case is bit more involved.  Clearly $n$ can be expressed as $n= (8m_1 + 4, 8k_1 + 4)$ for some $m_1, k_1 \in \mathbb{Z}$. Then
	$$
	3n =  6(4m_1 + 2, 4k_1 + 2).
	$$
	Let $\alpha_1 = 6$ and $\alpha_2 = (4m_1 + 2, 4k_1 + 2)$. That would imply (by Lemma \ref{lem2.1}) 
	\begin{equation}\label{eq5.2}
		a + 2r = (2m_1 + 4, 2k_1 + 1).
	\end{equation}
	In what follows we will apply Lemma \ref{lem3.1} (iv), (v), with $M, N \in \mathbb{Z}$.
	First, set $a = (12M + 4, 6N + 1)$, with $\text{Nm}(a)= 6$. Thus \eqref{eq5.2} implies that
	$$
	r = (m_1 - 6M, k_1 - 3N).
	$$
	Employing $ab + n = r^2$ and $d \equiv 10 \pmod{48}$ (see, (iii) of Lemma \ref{lem3.1}), we get
	\begin{align*}
		b = & \frac{1}{6}\left((m_1 - 6M)^2 + d(k_1 - 3N)^2 - 8m_1 - 4, 2(m_1 - 6M)(k_1 - 3N) - 8k_1 - 4\right) \times
		\\ & (12M + 4, -6N - 1).
	\end{align*}
	To ensure the existence of $b$ in $\mathbb{Z}[\sqrt{d}]$, we must have,
	$$
	(m_1, k_1) \equiv  (0, 0), (0, 1), (2, 0), (2, 2), (4, 1), (4, 2) \pmod{(6, 3)}.
	$$
%	Thus for $a = (12M + 4, 6N + 1)$, we have 
%	$$
%	(m_1, k_1) \equiv  (0, 0), (0, 1), (2, 0), (2, 2), (4, 1), (4, 2) \pmod{(6, 3)}.
%	$$
	As before, we assume $a = (12M + 4, 6N - 1)$, with $\text{Nm}(a)= 6$. Then we arrive at %using \eqref{eq5.2},
	%$$
	%r = (m - 6M, k - 3N + 1).
	%$$
	%As before utilizing $ab + n = r^2$,
	\begin{align*}
		b = &
		\dfrac{1}{6}\times((m_1 - 6M)^2 + d(k_1 - 3N + 1)^2 - 8m_1 - 4, 2(m_1 - 6M)(k_1 - 3N + 1) - 8k_1 - 4) \\ 
		& \times (12M + 4, -6N + 1).
	\end{align*}
	As $b\in\mathbb{Z}[\sqrt{d}]$, so that we have additional cases of $(m_1, k_1)$, where
	$$
	(m_1, k_1) \equiv (0, 2), (4, 0) \pmod{(6, 3)}.
	$$
%Similarly for $a = (12M + 4, 6N - 1)$: 
%$$
%	(m_1, k_1) \equiv (0, 2), (4, 0) \pmod{(6, 3)}.
%	$$

	Similarly, we set $a = (12M + 2, 6N + 1)$ with $\text{Nm}(a)=-6$ to get 
	\begin{align*}
		b =&  \dfrac{1}{-6}\times((m_1 + 1 - 6M)^2 + d(k_1 - 3N)^2 - 8m_1 -4, 2(m_1 + 1 - 6M)(k_1-3N) - 8k_1 - 4)\\
		&\times (12M + 2, -6N - 1).
	\end{align*}
	For $b$ to be in $\mathbb{Z}[\sqrt{d}]$, 
	$$
	(m_1, k_1) \equiv (1, 0), (1, 1), (3, 2), (5, 0),  (5, 2) \pmod {(6, 3)}.
	$$
%	For $a = (12M + 2, 6N + 1)$:
%	$$
%	(m_1, k_1) \equiv (1, 0), (1, 1), (3, 2), (5, 0),  (5, 2) \pmod {(6, 3)}.
%	$$
	Again we choose $a = (12M + 2, 6N - 1)$, with $\text{Nm}(a)=- 6$, which gives 
	\begin{align*}
		b = & \dfrac{1}{-6}\times((m_1 - 6M + 1)^2 + d(k_1 - 3N + 1)^2 - 8m_1 - 4, 2(m_1 - 6M + 1)(k_1 - 3N + 1) - 8k_1 - 4)\\ & \times (12M + 2, -6N + 1).
	\end{align*}
	Thus for $b \in \mathbb{Z}[\sqrt{d}]$, 
	$$
	(m_1, k_1) \equiv (1, 2), (3, 0), (3, 1) \pmod{(6, 3)}.
	$$
	
%Thus for $a = (12M + 2, 6N - 1)$,
%$$
%	(m_1, k_1) \equiv (1, 2), (3, 0), (3, 1) \pmod{(6, 3)}.
%	$$
	
Finally for $a=(12M-2, 6N-1)$ one gets the same values for $(m_1, k_1)$ as in the case  $a=(12M+2, 6N+1)$. %Finally, it turns out that this construction does not give a quadruple only for $(m_1, k_1)\equiv (2, 1)\pmod{(3,3)}$. 
This completes the proof of Theorem \ref{1.2}.
%	Analogously,  $a = (12M - 2, 6N + 1)$ with $\text{Nm}(a)= -6$ provides 
%	\begin{align*}
%		b =& \dfrac{1}{-6}\times((m + 3 - 6M)^2 + d(k - 3N)^2 - 8m - 4, 2(m + 3 - 6M)(k - 3N) - 8k - 4)\\
%		& \times (12M - 2, -6N - 1).
%	\end{align*}
%	We need to be restricted to
%	$$
%	(m, k) \equiv (3, 1) \pmod{(6, 3)}
%	$$
%	for $b \in \mathbb{Z}[\sqrt{d}]$. 
%	
%	Finally, we consider $a = (12M - 2, 6N - 1)$ satisfying $\text{Nm}(a)=- 6$. This gives us 
%	\begin{align*}
%		b =& \dfrac{1}{-6}\times((m + 3 - 6M)^2 + d(k - 3N + 1)^2 -8m - 4, 2(m + 3 - 6M)(k - 3N + 1) - 8k - 4)\\
%		& \times (12M - 2, -6N + 1).
%	\end{align*}
%	Now for $b \in \mathbb{Z}[\sqrt{d}]$, one gets  
%	$$
%	(m, k) \equiv (5, 0) \pmod{(6, 3)}.
%	$$
%	As dealt with in the previous cases, these values of $a, b$ and $ r$ will guarantee the existence of infinitely many $D(n)$-quadruples in $\mathbb{Z}[\sqrt{d}]$. 
% %This completes the proof.
 %\end{proof}
	
\section{Proof of Theorem \ref{1.3}}
The proof of Theorem \ref{1.3} goes along the lines of that of Theorem \ref{1.2}, except the factorization of $3n$. However, we provide the outlines of the proof for convenience to the readers. The notations $\alpha_1$ and $\alpha_2$ are as in \S\ref{s3}. Assume that $n = (4m + 2, 4k)$, where $m, k \in \mathbb{Z}$.
	
	\subsection*{Case I:  Both $m$ and $k$ are even} 
Let $M, N \in \mathbb{Z}$, and let 
	\begin{align}\label{eqf}
		3n &= 6(2m + 1, 2k) \nonumber \\ &=  (12M + 4, -6N - 1)(12M + 4, 6N + 1)(2m + 1, 2k) \textnormal{\hspace{0.5cm}(Using Lemma \ref{lem3.1}(iv))}\nonumber \\ & = \alpha_1\alpha_2,
	\end{align}
	where
	$$
	\begin{cases}
		\alpha_1 = (12M + 4, -6N - 1),\\
		\alpha_2 = (24Mm + 12M + 8m + 4 + d(12Nk + 2k), 24Mk + 8k + 12Nm + 2m + 6N + 1).
	\end{cases}
	$$
	Now, $a = 4(6a_1 + 1, 6b_1)$ with $a_1, b_1\in\mathbb{Z}$ and $\text{Nm}(a)=16$, which  gives 
	$$
	r = (6Mm + 6M + 2m + (d/2)(6Nk + k) - 12a_1, 6Mk + 2k + 3Nm + (m/2) - 12b_1)
	$$
	and
	\begin{align*}
		b =& \dfrac{1}{4}\Big\{(6Mm + 6M + 2m + (d/2)(6Nk + k) - 12a_1)^2 +d(6Mk + 2k + 3Nm + (m/2) -\\ &
		12b_1)^2 - 4m - 2, 2(6Mm + 6M + 2m + (d/2)(6Nk + k) - 12a_1)(6Mk + 2k + 3Nm +\\ &
		(m/2) - 12b_1) - 4k ) \times (6a_1 + 1, -6b_1)\Big\}.
	\end{align*}
	
	Now for $r, b \in \mathbb{Z}[\sqrt{d}]$, since $d \equiv 2 \pmod{4}$, we must have $m\equiv 2\pmod 4$. Assume that 
	$$
	(\alpha, \beta) =  (6Mm + 6M + 2m + (d/2)(6Nk + k) , 6Mk + 2k + 3Nm + m/2).
	$$
	Then $r=(\alpha-12a_1, \beta-12b_1)$.  
	
	Now if $a + 4b + 4r = 0$, then
	\begin{align*}
		4 + \alpha^2 + d\beta^2 - 4m - 2 + 4\alpha = 0, \\
		2\alpha\beta - 4k + 4\beta = 0.
	\end{align*}
	By Lemma \ref{lemm1}, we conclude that there exist only finitely many $\alpha$ and $\beta$ which satisfy the above system of equations. We now rewrite $\alpha$ and $\beta$ as follows, 
	\begin{align*}
		\alpha &= 6M(m + 1) + N(3dk) + 2m + (d/2)k \\
		\beta &=  6Mk + 3Nm + (m/2) + 2k.
	\end{align*}
	These can be written as
	\begin{equation*}
		\begin{pmatrix}
			\alpha - 2m - (d/2)k \\
			\beta - (m/2) - 2k 
		\end{pmatrix}
		=
		\begin{pmatrix}
			6(m + 1) & 3dk\\
			6k & 3m
		\end{pmatrix}
		\begin{pmatrix}
			M\\
			N
		\end{pmatrix}
		.
	\end{equation*}
	Since $m\equiv 2\pmod 4$, $k$ is even, and $d \equiv 2 \pmod{4}$, so that the determinant of 
	\begin{equation*}
		\begin{pmatrix}
			6(m + 1) & 3dk\\
			6k & 3m
		\end{pmatrix}
	\end{equation*}
	is non-zero. As we have infinitely many choices for $M$ and $N$, so that there exist infinitely many $\alpha$ and $\beta$ for which $a + 4b + 4r \neq 0$. Hence we can take such $M$ and $N$ for which $a + 4b + 4r \neq 0$. Using these values of $a, b$ and $r$,  we can get infinitely many quadruples with the property $D(n)$ from Lemma \ref{lem2.1}, since we have infinitely many choices of $a$, by using Lemma \ref{lem3.1} $(i)$ and for checking the condition of non-zero and distinct elements of the set $\{a, b, a+ b + 2r, a+ 4b + 4r\}$ (given in Lemma \ref{lem2.1}), we use Lemma \ref{lemm1}. 
	
	In the case $m \equiv  0 \pmod{4}$, we replace  $n$ by $n = (16m_1 + 2, 8k_1)$ and then consider \eqref{eqf} with
	\begin{align*}
		\alpha_1 = & (-12M - 2, 6N + 1),\\
		\alpha_2 = & (96Mm_1 + 12M +16m_1 + 2 + d(24Nk_1 + 4k_1), 48Mk_1 + 8k_1 + 48Nm_1 + 8m_1 + \\
		 &  6N + 1),
	\end{align*}
where $m_1, k_1 \in \mathbb{Z}$.
	This gives by utilizing $a = (12a_1 + 4, 6b_1 + 1)$ with $a_1, b_1\in\mathbb{Z}$ and $\text{Nm}(a)=6$,
	$$
	r = (24Mm_1 + 4m_1 + d(6Nk_1 + k_1) - 6a_1 - 2, 12Mk_1 + 2k_1 + 12Nm_1 + 2m_1 + 3N - 3b_1)
	$$
	and 
	\begin{align*}
		b =& \dfrac{1}{6}\Big\{((24Mm_1 + 4m_1 + d(6Nk_1 + k_1) - 6a_1 - 2)^2 + d(12Mk_1 + 2k_1 + 12Nm_1 + \\ & 2m_1 + 3N -  3b_1)^2
		- 16m_1 - 2,  2(24Mm_1 + 4m_1 + d(6Nk_1 + k_1) - 6a_1 - 2)(12Mk_1 + \\ & 2k_1 +  12Nm_1 + 2m_1 
		+ 3N - 3b_1) - 8k_1) \times (12a_1 + 4, -6b_1 - 1)\Big\}.
	\end{align*}
	Using $d \equiv 10 \pmod{48}$ (from Lemma \ref{lem3.1}(iii)), these further imply that
	$$
	(m_1, k_1) \equiv (0, 1), (0, 2), (1, 0), (1, 1), (2, 0), (2, 2) \pmod{(3, 3)}.
	$$
	
	Similarly, for $a = (12a_1 - 4, 6b_1 + 1)$ with $a_1, b_1\in\mathbb{Z}$ and $\text{Nm}(a)=6$, we have
	$$
	r = (24Mm_1 +4m_1 + d(6Nk_1 + k_1) - 6a_1 + 2, 12Mk_1 + 2k_1 + 12Nm_1 + 2m_1 + 3N - 3b_1)
	$$
	and
	\begin{align*}
		b =& \dfrac{1}{6}\Big\{((24Mm_1 +4m_1 + d(6Nk_1 + k_1) - 6a_1 + 2)^2 + d(12Mk_1 + 2k_1 + 12Nm_1 + 2m_1 + \\ & 3N - 3b_1)^2  - 16m_1 - 2, 2(24Mm_1 +4m_1 + d(6Nk_1 + k_1) - 6a_1 + 2)(12Mk_1 + \\ & 2k_1 + 12Nm_1 +  2m_1 + 3N - 3b_1) - 8k_1) \Big\} \times \Big(12a_1 - 4, - 6b_1 - 1\Big).
	\end{align*}
	For  $b$ to be in $\mathbb{Z}[\sqrt{d}]$,
	$$
	(m_1, k_1) \equiv (1, 2) \pmod{(3, 3)}.
	$$
	
	The factorization \eqref{eqf} with 
	$$
	\begin{cases}
		\alpha_1 = & (12M + 2, 6N + 1),\\
		\alpha_2 = & (-96Mm_1 - 12M -16m_1 - 2 +d(24Nk_1 + 4k_1), -48Mk_1 - 8k_1 + 48m_1N + \\ & 8m_1 + 6N + 1),
	\end{cases} 
	$$
	as well as $a = (12a_1 + 4, 6b_1 - 1)$ with $a_1, b_1\in\mathbb{Z}$ and $\text{Nm}(a)= 6$ provides
	$$
	r = (-24Mm_1 -4m_1 +d(6Nk_1 + k_1) - 6a_1 - 2, -12Mk_1 - 2k_1 + 12m_1N + 2m_1 + 3N + 1 - 3b_1)
	$$
	and
	\begin{align*}
		b =& \dfrac{1}{6}\Big\{(-24Mm_1 -4m_1 +d(6Nk_1 + k_1) - 6a_1 - 2)^2 + d(-12Mk_1 - 2k_1 + 12m_1N + 2m_1 + 3N + 1 \\ &- 3b_1)^2 - 16m_1 - 2, 2(-24Mm_1 -4m_1 +d(6Nk_1 + k_1) - 6a_1 - 2)(-12Mk_1 - 2k_1 + 12m_1N \\ & + 2m_1 + 3N + 1 - 3b_1) - 8k_1\Big\}\times\Big(12a_1 + 4, -6b_1 + 1\Big).
	\end{align*}

For $b \in \mathbb{Z}[\sqrt{d}]$,
	$$
	(m_1, k_1) \equiv (2, 1) \pmod{(3, 3)}.
	$$ 
	Finally, owing to Lemma \ref{lem3.1}, there are infinitely many choices of $M$ and $N$, and hence there are infinitely many choices for such $a, b$ and $r$.
	
To conclude this case,  we have covered all possibilities for $(m_1, k_1)$, except $(m_1, k_1) \not\equiv (0, 0) \pmod{(3, 3)}$. Hence, there exist infinitely many Diophantine quadruples in $\mathbb{Z}[\sqrt{d}]$ with the property $D(16m_1 + 2, 8k_1)$, where $(m_1, k_1) \not\equiv (0, 0) \pmod{(3, 3)}$.
	
	\subsection*{Case II: $m$ is even and $k$ is odd} In this case too we work with the factorization \eqref{eqf}.  We use
	$$
	\begin{cases}
		\alpha_1 = (12M + 4, -6N - 1),\\
		\alpha_2 = (24Mm + 12M + 8m + 4 + d(12Nk + 2k), 24Mk + 8k + 12Nm + 2m + 6N + 1)
	\end{cases}
	$$
	
	and $a = 2(6a_1 + 1, 6b_1)$ with $a_1, b_1\in\mathbb{Z}$ and $\text{Nm}(a)=4$. These provide us, 
	$$
	r =  (6Mm + 6M + 2m + 2 + (d/2)(6Nk + k) - 6a_1 - 1, 6Mk + 2k + 3Nm + (m/2) - 6b_1)
	$$
	and
	\begin{align*}
		b =& \dfrac{1}{2}\Big\{(6Mm + 6M + 2m + 2 + (d/2)(6Nk + k) - 6a_1 - 1)^2 + d(6Mk + 2k + 3Nm +\\& (m/2) - 6b_1)^2 - 4m - 2, 2(6Mm + 6M + 2m + 2 + (d/2)(6Nk + k) - 6a_1 - 1)(6Mk +\\ & 2k + 3Nm + (m/2) - 6b_1) - 4k\Big\}\times \Big(6a_1 + 1, -6b_1\Big).
	\end{align*}
	
	\subsection*{Case III: $m$ is odd and $k$ is even} Here, we use \eqref{eqf} with 
	\begin{align*}
		\alpha_1& = (-12M - 2, 6N + 1),\\
		\alpha_2& = (24Mm + 12M +4m + 2 + d(12Nk + 2k), 24Mk + 4k + 12Nm + 2m + 6N + 1).
	\end{align*}
	Then, $a = 2(6a_1 + 1, 6b_1)$ with $a_1, b_1\in\mathbb{Z}$ and $\text{Nm}(a)=4$ gives
	$$
	r =  (12Mm +2m + d(6Nk + k), 12Mk + 2k + 6Nm + m + 6N + 1)
	$$
	and
	\begin{align*}
	b = & \dfrac{1}{2}\Big\{(12Mm +2m + d(6Nk + k))^2 + d(12Mk + 2k + 6Nm + m + 6N + 1)^2 - 4m -\\ & 2, 2(12Mm +2m + d(6Nk + k))(12Mk + 2k + 6Nm + m + 6N + 1) - 4k\Big\} \times \\ & \Big(6a_1 + 1, -6b_1\Big).
	\end{align*}
	
	\subsection*{Case IV: Both $m$ and $k$ are odd} 
	
	The choices of $\alpha_1$ and $\alpha_2$ as in Case III work in this case too. 
	We set $a = 4(6a_1 + 1, 6b_1)$ with  $a_1, b_1\in\mathbb{Z}$ and $\text{Nm}(a)=16$ to get  
	$$
	r = (6Mm + m + (d/2)(6Nk + k) - 12a_1 - 2, 6Mk + k + 3Nm + (m + 1)/2 + 3N - 12b_1)
	$$
	and
	\begin{align*}
		b = & \dfrac{1}{4}\Big\{((6Mm + m + (d/2)(6Nk + k) - 12a_1 - 2)^2 + d(6Mk + k + 3Nm + (m + 1)/2 + \\ &
		3N - 12b_1)^2 - 4m - 2, 2(6Mm + m + (d/2)(6Nk + k) - 12a_1 - 2)(6Mk + k + 3Nm\\
		&
		+ (m + 1)/2 + 3N - 12b_1) -4k) (6a_1 + 1, -6b_1)\Big\}.
	\end{align*}
	These would imply $m\equiv 3\pmod 4$ whenever $r, b \in \mathbb{Z}[\sqrt{d}]$. The existence of infinitely many quadruples
	can be seen by similar argument of $n = (4m + 2, 4k)$ in Case I with $m \equiv 2 \pmod4$ and even $k$.

	The next case is $m\equiv 1\pmod 4$ and here $n$ can be replaced by  $n= (16m_1 + 6, 8k_1 + 4)$ with $m_1, k_1 \in \mathbb{Z}$. The factorization uses in this case is:
	\begin{equation}\label{eq3}
		3n = \alpha_1\alpha_2,
	\end{equation}
	where,
	\begin{align*} 
		\alpha_1 =& (12M + 4, -6N - 1), \\
		\alpha_2 =
		& (96Mm_1 + 36M + 32m_1 + 12 + d(24Nk_1 + 12N + 4k_1 + 2), 48Mk_1 + 24M + 16k_1 \\
		&+ 11 + 48Nm_1 + 18N + 8m_1).
	\end{align*}
	We set $a = (12a_1 + 2, 6b_1 + 1)$ with $a_1, b_1\in\mathbb{Z}$ and $\text{Nm}(a)= -6$, which gives
	\begin{align*}
			r = & (24Mm_1 + 12M + 8m_1 + 4 + (d/2)(12Nk_1 + 6N + 2k_1 + 1) - 6a_1 - 1, 12Mk_1 + \\ & 6M + 4k_1 - 3a_1 + 2 + 12Nm_1 + 3N + 2m_1)
	\end{align*}
	 and
	 \begin{align*}
	 	b = & \dfrac{1}{-6}\Big\{(24Mm_1 + 12M + 8m_1 + 4 + (d/2)(12Nk_1 + 6N + 2k_1 + 1) - 6a_1 - 1)^2 + \\ & d(12Mk_1 + 6M + 4k_1 - 3a_1 +   2 + 12Nm_1 + 3N + 2m_1)^2 - 16m_1 - 6, 2(24Mm_1 + \\ & 12M + 8m_1 + 4 + (d/2)(12Nk_1 + 6N + 2k_1 + 1) - 6a_1 - 1)(12Mk_1 + 6M + 4k_1 - \\ & 3a_1 +  2 + 12Nm_1 + 3N + 2m_1) - 8k_1 - 4\Big\} \times \Big(12a_1 + 2, -6b_1 - 1\Big).
	 \end{align*}
Now, using $d \equiv 10 \pmod{48}$ (from Lemma \ref{lem3.1}(iii)),
	$$
	(m_1, k_1) \equiv (0, 0), (0, 1), (1, 1), (2, 0), (2, 2) \pmod{(3, 3)},
	$$
	for $b \in \mathbb{Z}[\sqrt{d}]$.

	Similarly, $a = (12a_1 -2, 6b_1 + 1)$ with $a_1, b_1\in\mathbb{Z}$ and $\text{Nm}(a)= -6$ provides
	\begin{align*}
		r = & (24Mm_1 + 12M + 8m_1 + 5 + (d/2)(12Nk_1 + 6N + 2k_1 + 1) - 6a_1, 12Mk_1 + 6M + \\ & 4k_1 + 2 + 12Nm_1 + 3N + 2m_1 - 3b_1)
	\end{align*}
	and
	\begin{align*}
		b = & \dfrac{1}{-6}\Big\{(24Mm_1 + 12M + 8m_1 + 5 + (d/2)(12Nk_1 + 6N + 2k_1 + 1) - 6a_1)^2 + \\ & d(12Mk_1 +  6M +  4k_1 + 2 + 12Nm_1 + 3N + 2m_1 - 3b_1)^2 - 16m_1 - 6, 2(24Mm_1 + \\ & 12M + 8m_1 + 5 +  (d/2)(12Nk_1 + 6N + 2k_1 + 1) - 6a_1)(12Mk_1 + 6M + 4k_1 + 2 + \\ & 12Nm_1 + 3N + 2m_1 -  3b_1) - 8k_1 - 4\Big\} \times \Big(12a_1 - 2, -6b_1 - 1\Big).
	\end{align*}
	For $b \in \mathbb{Z}[\sqrt{d}]$,
	  $$(m_1, k_1) \equiv (0, 2) \pmod{(3, 3)}.$$
	
	Again, we use \eqref{eq3} by taking 
	\begin{align*}
		\alpha_ 1 =& (12M + 4, 6N +1), \\ 
		\alpha_2 =& (96Mm_1 + 36M + 32m_1 + 12 + d(-24Nk_1 - 12N - 4k_1 - 2), 48Mk_1 + 24M + \\ & 16k_1 + 8 -
		48Nm_1 -18N - 8m_1 - 3).
	\end{align*}
	Then we choose 
	$a = (12a_1 + 2, 6b_1 - 1)$ with $a_1, b_1\in\mathbb{Z}$ and $\text{Nm}(a)=-6$, which gives
	\begin{align*}
	r = & (24Mm_1 + 12M + 8m_1 + 3 + (d/2)(-12Nk_1 - 6N - 2k_1 - 1) - 6a_1, 12Mk_1 + 6M + 4k_1 + 2 -\\ &
	12Nm_1 -3N - 2m_1 - 3b_1)
	\end{align*}
and
\begin{align*}
	b = & \dfrac{1}{-6}\Big\{(24Mm_1 + 12M + 8m_1 + 3 + (d/2)(-12Nk_1 - 6N - 2k_1 - 1) - 6a_1)^2 + d(12Mk_1 + \\ & 6M + 4k_1 + 2 -
	12Nm_1 -3N - 2m_1 - 3b_1)^2 - 16m_1 - 6, 2(24Mm_1 + 12M + 8m_1 + 3 + \\ & (d/2)(-12Nk_1 - 6N - 2k_1 - 1) - 6a_1)(12Mk_1 + 6M + 4k_1 + 2 - 12Nm_1 -3N - 2m_1 - \\ &  3b_1) - 8k_1 - 4\Big\} \times \Big(12a_1 + 2, -6b_1 + 1\Big).
\end{align*}
 For $b \in \mathbb{Z}[\sqrt{d}]$, 
$$(m_1, k_1) \equiv (1,0) \pmod{(3, 3)}.$$
The existence of infinitely many $D(n)$-quadruples in $\mathbb{Z}[\sqrt{d}]$ is guaranteed by the above choices of $a, b$ and $r$ in each case. 

	To conclude this case,  we have covered all possibilities for $(m_1, k_1)$ except $(m_1, k_1) \not\equiv (2, 1) \pmod{(3, 3)}$. Therefore, there exist infinitely many Diophantine quadruples in $\mathbb{Z}[\sqrt{d}]$ with the property $D(16m_1 + 6, 8k_1 + 4)$, where $(m_1, k_1) \not\equiv (2, 1) \pmod{(3, 3)}$.
	
\section{Concluding Remarks}
Given a square-free integer $d\equiv 2\pmod 4$, the existence of  $D(n)$-quadruples in the ring $\mathbb{Z}[\sqrt{d}]$ for some $n\in \mathbb{Z}[\sqrt{d}]$ has been investigated in \cite{CGH2022, FR2004}. We investigate this problem for the remaining values of $n$. However, our method does not work for a few values of $n$, i.e., $ n \in \{4(12r+5, 6s+3), 4(12r+11, 6s+3), (48r+38,24s+12), (48r+2, 24s)\}$ with $r, s\in \mathbb{Z}$.% Indeed, for showing the existence of quadruples with the property $D(n)$, we have used Lemma \ref{lem2.1}. According to this lemma,  $ab + n = r^2$ and $3n = \alpha_1\alpha_2$. This forces that we have two equations, namely,  $r = \dfrac{\alpha_1 + \alpha_2}{4} - \dfrac{a}{2}$ and $ab + n = r^2$. 
	%This imples that $a + 2r = \dfrac{\alpha_1 + \alpha_2}{2}$. 
	%If one can handle these equations, then we can find the existence of Diophantine quadruples with the property $D(n)$. Now, for ceratin $d$'s, we discuss Diophantine quadruples  in $\mathbb{Z}[\sqrt{d}]$ with the property $D(n)$, with such $n$'s which are not considered in Theorems \ref{1.2} and \ref{1.3}.   

	%The authors believe that for these exception cases of $n$, there does not exist any Diophantine quadruple in $\mathbb{Z}[\sqrt{d}]$ with the property $D(n)$. So this can be give another family Diophantine quadruples in $\mathbb{Z}[\sqrt{d}]$ with the property $D(n)$, for which Conjecture \ref{Con1.1} does not hold. Unfortunately, we have no example of Diophantine quadruple in $\mathbb{Z}[\sqrt{d}]$ with the property $D(n)$, for these $n$'s.
	
	We discuss some examples for the existence of $D(n)$-quadruples in $\mathbb{Z}[\sqrt{d}]$ for  these exceptions. We first shorten these exceptions with the help of  \cite[Theorem 1.1]{CGH2022}, and then we provide some examples for the remaining cases.
 
	Let $d = 2N$ such that \eqref{eqi1} and \eqref{eqi2} are solvable in integers, where  $N\in \mathbb{N}$. Assume that $n = 4(12m + 5, 6k+ 3)$ with $m = \alpha N + \beta$ and $k = \alpha_1 N + \beta_1$, where $\alpha, \beta, \alpha_1, \beta_1 \in \mathbb{Z}$. Then $ n =  4(12 \alpha N + 12 \beta + 5, 6\alpha_1 N + 6 \beta_1 + 3)$. Utilizing (iii) of Lemma \ref{lem3.1}, we get $2, 3 \nmid N$ and thus we can choose $\beta, \beta_1$ such that $12\beta + 5$ and $6\beta_1 + 3$ are of the form $N\gamma$ and $N\gamma_1$,  respectively with odd integers $\gamma$ and $\gamma_1$. Thus $n =  2N(24\alpha + 2\gamma, 12\alpha_1 + 2\gamma_1)$, since $2\gamma, 2\gamma_1 \equiv 2 \pmod 4$, so that $24 \alpha + 2\gamma$ and  $12\alpha_1 + 2\gamma_1$ are of the form $4t_1+2$ for some integer $t_1\geq 1$. 
 
 Again $2N$ is square in $\mathbb{Z}[\sqrt{d}]$, and thus \cite[Theorem 1.1]{CGH2022} and  Lemma \ref{lemp} together show that there exist infinitely many $D(n)$-quadruples in $\mathbb{Z}[\sqrt{d}]$. Analogously, we can draw a similar conclusion for $n = 4(12m + 11, 6k + 3)$.
	We now consider $n = (4(12m + 9) + 2, 4(6k + 3))$. As in the above, $n = 2(24N\alpha + 24\beta + 19, 12N\alpha_1 + 12\beta_1 + 6)$. Since $2, 3 \nmid N$, so that we can choose $\beta, \beta_1$ such that $24\beta + 19$ and $12\beta_1 + 6$ are of the form $N\gamma$ and $N\gamma_1$, respectively. Using (iii) of Lemma \ref{lem3.1}, we get $N \equiv 1 \pmod 4$, and thus  $\gamma \equiv 3 \pmod{4}$ and $\gamma_1 \equiv 2 \pmod 4$. Finally we use \cite[Theorem 1.1]{CGH2022} and  Lemma \ref{lemp} to conclude that there exist infinitely many $D(n)$-quadruples in $\mathbb{Z}[\sqrt{d}]$. Analogously we can establish the same for $n = (4(12m) + 2, 4(6k))$.
	%and $(4(12m) + 2, 4(6k))$.
	%From Lemma \ref {lem3.1} $(iii)$, $p$ must be an odd integer. .............
	
	We now provide some examples supporting  the existence of  $D(n)$-quadruple in $\mathbb{Z}[\sqrt{10}]$ for the exceptional values of $n$. 
	\subsection*{Example 1}
	We consider $d = 10$  and $n = 4(12m + 5, 6k + 3)$ with $m, k \in \mathbb{Z}$. Let $m = 5M$ and $k = 5K + 2$, where $M, K \in \mathbb{Z}$. Then $n = 4(5(12M + 1), 30K + 15)$, which can be written as  $n= 10(24M + 2, 12K + 6)$. Thus $n$ is of the form $10(4m' + 2, 4k' + 2)$ with $m', k' \in \mathbb{Z}$. Therefore using \cite[Theorem 1.1]{CGH2022} and  Lemma \ref{lemp}, we conclude that there exist infinitely many $D(n)$-quadruples in $\mathbb{Z}[\sqrt{10}]$. Analogously, we can show the same for  $n = 4(12m + 11, 6k + 3)$ by putting $m = 5M + 2$ and $k = 5K + 2$. 
	
	\subsection*{Example 2}
	Suppose $d = 10$ and $n = (4(12m + 9) + 2, 4(6k + 3)) = 2(24m + 19, 12k + 6)$. Let $m = 5M + 4$ and $k = 5K + 2$. Then $n = 10(24M + 23, 12k + 6)$. Since $24M + 23 \equiv 3 \pmod{4}$ and $12K + 6 \equiv 2 \pmod 4$, so that by \cite[Theorem 1.1]{CGH2022} and  Lemma \ref{lemp}, we can conclude that there exist infinitely many  $D(n)$-quadruples in $\mathbb{Z}[\sqrt{10}]$. Similar conclusion can be drawn for $n = (4(12m) + 2, 4(6k))$ by taking $m \equiv 1 \pmod 5$ and $k \equiv 0 \pmod 5$.
	
	\subsection*{Example 3}
	Assume that $d = 10$ and $n = 4(12m + 5, 6k + 3)$ with $m \equiv  2, 3 \pmod 5$.
	We factorize $3n$ as follows:
	\begin{align*}
		3n &= 12(12m + 5, 6k + 3) \\
		&= (-18, 6)(3, 1)(24m + 10, 12k + 6)\\
		&= (-18, 6)(120k + 72m + 90, ~36k + 24m + 28).
	\end{align*}
We take $\alpha_1$ and $\alpha_2$ to be the first and the second factor of the above equation, respectively. Further utilizing Lemma \ref{lem2.1} we get
	$$
	a + 2r = (60k + 36m + 36, 18k + 12m + 17).
	$$
	We choose $a = (19, 6)^t(0, 1)$ with $\text{Nm}(a)=-10$, where $t \in \mathbb{N}$. This implies that  there exist $\alpha, \beta  \in \mathbb{Z}$ such that $a = (20\alpha, 10\beta - 1)$, and thus  $ r = (30k + 18m - 10\alpha + 18, ~9k + 6m + 9 - 5\beta)$. Further $ab + n = r^2$ implies
	$$
	b = \dfrac{(r^2 - n)(20\alpha, -10\beta + 1)}{-10}.
	$$
	Since $m \equiv 2,\text{~~or~~} 3 \pmod 5$, $b \in \mathbb{Z}[\sqrt{10}]$ and  we have infinitely many $a$'s, therefore by using Lemmas \ref{lem2.1} and \ref{lemm1}, we get  infinitely many  $D(n)$-quadruples in $\mathbb{Z}[\sqrt{10}]$.
	Analogously, we can show the existence of infinitely many   $D(n)$-quadruples in $\mathbb{Z}[\sqrt{10}]$ for $n = 4(12m + 11, 6k + 3)$ when $m \equiv 0,\text{~~or~~} 4 \pmod 5$.
	
	\subsection*{Example 4}
	Suppose that $d = 10$ and $n = (4(12m + 9) + 2,~4(6k + 3))$ with $m \equiv 1,\text{~~or~~} 2 \pmod 5$. We factorize 
	\begin{align*}
		3n &= (4, 1)(4, -1)(2(12m + 9) + 1,~2(6k + 3)) \\
		&= (4, 1)(-120k + 96m + 16,~48k - 24m + 5).
	\end{align*}
	We choose $\alpha_1$ and $\alpha_2$ to be the first and the second factor of the last equation, respectively. We use Lemma \ref{lem2.1} to get  $a + 2r = (-60k + 48m + 10,~24k - 12m + 3)$. Let $a = (19, 6)^t(10, 3)$ with $\text{Nm}(a)=10$, where $t \in \mathbb{N}$. Thus there exist $\alpha,~\beta \in \mathbb{Z}$ such that $a = (20\alpha + 10, 10\beta + 3)$ and thus $r = (24m - 30k - 10\alpha, -6m + 12k - 5\beta)$. Therefore Lemma \ref{lem2.1} gives
	$$
	b = \dfrac{r^2 - n}{a}.
	$$
	Since $m \equiv 1,\text{~~or~~} 2 \pmod 5$, so that $b \in  \mathbb{Z}[\sqrt{10}]$. Hence there exist infinitely many   $D(n)$-quadruples in $\mathbb{Z}[\sqrt{10}]$. Analogously, we can construct $D(n)$-quadruples for $n = (4(12m) + 2, 4(6k))$ when $m \equiv 3,\text{~~or~~} 4 \pmod 5$.
	
	The problem of existence of infinitely many $D(n)$-quadruples in $\mathbb{Z}[\sqrt{10}]$ for $n\in \mathbb{Z}[\sqrt{10}]$ is solved, except for $n\in \mathcal{S}_0:=S_1\cup S_2\cup S_3\cup S_4$, where  \\ $S_1=\{4(12m+5, 6k+3): (m, k)\equiv (0,0), (0,1), (0,3), (0,4)\pmod{(5,5)} \text{ or } m\equiv 1,4\pmod 5\}$,\\
	$S_2=\{4(12m+11, 6k+3): (m, k)\equiv (2,0), (2,1), (2,3), (2,4)\pmod{(5,5)} \text{ or } m\equiv 1,3\pmod 5\}$,\\
	$S_3=\{(4(12m+9)+2, 4(6k+3)): (m, k)\equiv (4,0), (4,1), (4,3), (4,4)\pmod{(5,5)} \text{ or } m\equiv 0,3\pmod 5\}$ and\\
	$S_4=\{48m+2, 36k): (m, k)\equiv (1,1), (1,2), (1,3), (1,4)\pmod{(5,5)} \text{ or } m\equiv 0,2\pmod 5\}$.
	
	Finally, we put the following question for $n\in \mathcal{S}_0$. 
	\begin{q}
		Do there exist infinitely many $D(n)$-quadruples in $\mathbb{Z}[\sqrt{10}]$ when $n \in S_0$?
		%From Example 1, 2, 3, and 4, one can see that still some $n$'s are remaining for analyzing of  Diophantine quadruples in $\mathbb{Z}[\sqrt{10}]$ with the property $D(n)$. Unfortunately, we have no example for quadruples with the property $D(n)$ for these $n$'s. From here, there may be exist a different family of Diophantine quadruples in $\mathbb{Z}[\sqrt{10}]$ from \cite{CGH2022}, for which Conjecture \ref{Con1.1} does not hold.
	\end{q}

		\section*{Acknowledgements}
		The authors would like to thank the anonymous referees for their valuable suggestions/comments that immensely improved the presentation of the paper.  A. Hoque acknowledges SERB MATRICS Project (MTR/2021/000762) and SERB CRG Project CRG/2023/007323, Govt. of India.


\begin{thebibliography}{99}
			\bibitem{MR2004} F. S. Abu Muriefah and A. Al Rashed, \textit{Some Diophantine quadruples in the ring $\mathbb{Z}[\sqrt{-2}]$}, Math. Commun. \textbf{9} (2004), 1--8.
			
			\bibitem{A2019} N. Ad\v{z}aga, {\it On the size of Diophantine $m$-tuples in imaginary quadratic number rings}, Bull. Math. Sci. {\bf 9} (2019), no. 3, Article ID: 1950020, 10pp.
			\bibitem{BD1969} A. Baker and H. Davenport, {\it The equations $3x^2 - 2 = y^2$ and $8x^2 - 7 = z^2$}, Quart. J. Math. Oxford Ser. (2) {\bf 20} (1969), 129--137.
			\bibitem{BTF2019} M. Bliznac Trebje\v{s}anin and A. Filipin, {\it Nonexistence of $D(4)$-quintuples}, J. Number Theory {\bf 194} (2019), 170--217.
			
			\bibitem{BCM2020} N. C. Bonciocat, M. Cipu and M. Mignotte, {\it There is no Diophantine $D(-1)$-quadruple}, J. London Math. Soc. {\bf 105} (2022), 63--99.
			\bibitem{BR1985} E. Brown, {\it Sets in which $xy + k$ is always a square}, Math. Comp. {\bf 45} (1985), 613--620.
			
			
		
			
			\bibitem{CGH2022} K. Chakraborty, S. Gupta, and A. Hoque, \textit{On a conjecture of Franu\v si\'c and Jadrijevi\' c: Counter-examples}, Results Math. {\bf 78} (2023), no. 1, 14pp, article no. 18.
			
			\bibitem{CGH22} K. Chakraborty, S. Gupta and A. Hoque, {\it Diophantine triples with the property $D(n)$ for distinct $n$'s}, Mediterr. J. Math. {\bf 20} (2023), no. 1, 13pp, article no. 31.
			
			
			\bibitem{DU1993} A. Dujella, {\it Generalization of a problem of Diophantus}, Acta Arith. \textbf{65} (1993), 15--27.
			
			\bibitem{DUJE1996} A. Dujella, \textit{Some polynomial formulas for Diophantine quadruples}, Grazer Math. Ber. \textbf{328} (1996), 25–30.
			
			\bibitem{DU1997} A. Dujella, {\it The problem of Diophantus and Davenport for Gaussian integers}, Glas. Mat. Ser. III {\bf 32} (1997), 1--10.
			
			\bibitem{DU2004} A. Dujella, {\it There are only finitely many Diophantine quintuples}, J. Reine Angew. Math. \textbf{566} (2004), 183--214.
			
			\bibitem{DU21}  A. Dujella, \textit{Number Theory}, \v{S}kolska knjiga, Zagreb, 2021.
						\bibitem{DU2024} A. Dujella, \textit{Diophantine $m$-tuples and Elliptic Curves},  Springer, Cham, 2024.

			\bibitem{DU23} A. Dujella, {\it Triples, quadruples and quintuples which are $D(n)$-sets for several $n$'s},  in: Class Groups of Number Fields and Related Topics, K. Chakraborty, A. Hoque and P. P. Pandey (eds.) (to appear). 
			
			
			
			
			\bibitem{EFF2014} C. Elsholz, A. Filipin and Y. Fujita, {\it On Diophantine quintuples and $D(-1)$-quadruples}, Monats. Math. {\bf 175} (2014), 227--239.
			
				\bibitem{FR2004} Z. Franu\v si\'c, \textit{Diophantine quadruples in the ring $\mathbb{Z}[\sqrt{2}]$}, Math. Commun. \textbf{9} (2004), 141--148. 
			
			%\bibitem{DF2007} A. Dujella and Franu\v si\'c, \textit{On differences of two Squares in some quadratic fields},  Rocky Mountain J. Math. \textbf{37} (2007), 429--453. 
			
			\bibitem{FR2008} Z. Franu\v si\'c, \textit{Diophantine quadruples in $\mathbb{Z}[\sqrt{4k + 3}]$}, Ramanujan J. \textbf{17} (2008), 77--88.
			
			\bibitem{FR2009} Z. Franu\v si\'c, \textit{A Diophantine problem in $\mathbb{Z}[\sqrt{(1 + d)/2}]$}, Studia Sci. Math. Hungar. \textbf{46} (2009), 103--112.
			
			\bibitem{FR2013} Z. Franu\v si\'c, \textit{Diophantine quadruples in the ring of integers of the pure cubic field $\mathbb{Q}(\sqrt[3]{2})$}, Miskolc Math.
			Notes \textbf{14} (2013), 893--903.
			\bibitem{FS2014} Z. Franu\v si\'c and I. Soldo, \textit{The problem of Diophantus for integers of $\mathbb{Q}(\sqrt{-3})$}, Rad Hrvat. Akad. Znan. Umjet. Mat. Znan. \textbf{18} (2014), 15--25.
			
			\bibitem{FJ2019} Z. Franu\v si\'c and B. Jadrijevi\' c, \textit{$D(n)$-quadruples in the ring of integers of $\mathbb{Q}(\sqrt{2}, \sqrt{3})$}, Math. Slovaca \textbf{69} (2019), 1263--1278.
			
			
			\bibitem{G2021} S. Gupta, {\it $D(-1)$ tuples in imaginary quadratic fields}, Acta Math. Hungar. {\bf 164} (2021), 556--569.
			
			%\bibitem{HT2011} B. He and A. Togb\'e, {\it On the $D(-1)$-triple $\{1, k^2+1, k^2+2k + 2\}$ and its unique $D(1)$-extension}, J. Number Theory {\bf 131} (2011), 120--137.
			
			\bibitem{HTZ2019} B. He, A. Togb\'e and V. Ziegler, {\it There is no Diophantine quintuple}, Trans. Amer. Math. Soc. \textbf{371} (2019), 6665--6709.
			
			
			
		
			
			
			
			\bibitem{MA2012} Lj. Juki\'c Mati\'c, {\it Non-existence of certain Diophantine quadruples in rings of integers of pure cubic fields}, Proc. Japan Acad. Ser. A  Math. Sci. {\bf 88} (2012), no. 10, 163--167.
			
			%\bibitem{MR1985} S. P. Mohanty and M. S. Ramamsamy, {\it On $P_{r,k}$ sequences}, Fibonacci Quart. {\bf 23} (1985), 36--44.
			
			
			\bibitem{SO2013} I.  Soldo, \textit{On the existence of Diophantine quadruples in $\mathbb{Z}[\sqrt{-2}]$}, Miskolc Math. Notes \textbf{14} (2013), 265--277.
			
				
		\end{thebibliography}
	\end{document}